\documentclass[letterpaper,10pt,journal,twoside]{IEEEtran}

\usepackage{graphicx} % for pdf, bitmapped graphics files
\usepackage{epsfig} % for postscript graphics files
\usepackage{amsmath} % assumes amsmath package installed
\usepackage{amssymb}  % assumes amsmath package installed
\usepackage{mathtools}
\usepackage{tabularx}
\allowdisplaybreaks
\usepackage{multirow}

\usepackage{epstopdf}

\usepackage{microtype} %typographi
\usepackage{float}
\newtheorem{remark}{Remark}
\newtheorem{theorem}{Theorem}

\newtheorem{assumption}{Assumption}

\newtheorem{lemma}{Lemma}
\newtheorem{definition}{Definition}

\usepackage{xspace}
\usepackage{soul}
\usepackage{mathtools}
\usepackage{makecell}
\usepackage{cite}
\usepackage{todonotes}
\usepackage{nicefrac}
\usepackage{fancyhdr}

\newcommand{\ra}{\rightarrow}

\newcommand{\ie}{\unskip, i.\,e.,\xspace}
\newcommand{\eg}{\unskip, e.\,g.,\xspace}

\newcommand{\N}{\ensuremath{\mathbb{N}}}

\newcommand{\R}{\ensuremath{\mathbb{R}}}

\newcommand{\X}{\ensuremath{\mathbb{X}}}
\newcommand{\Y}{\ensuremath{\mathbb{Y}}}
\newcommand{\F}{\ensuremath{\mathbb{F}}}
\newcommand{\U}{\ensuremath{\mathbb{U}}}

\newcommand*\diff{\mathop{}\!\mathrm{d}}
\newcommand{\eps}{\ensuremath{\varepsilon}}

\newcommand{\spc}{\ensuremath{\,\,}}

\newcommand{\ball}{\ensuremath{\mathcal B}}
\newcommand*\circled[1]{\tikz[baseline=(char.base)]{\node[shape=circle,draw,inner sep=1pt] (char) {#1};}}

\newcommand{\co}{\ensuremath{\overline{\text{co}}}}

\begin{document}
\title{Practical sample-and-hold stabilization of nonlinear systems under approximate optimizers}

\author{Pavel Osinenko, Lukas Beckenbach and Stefan Streif% <-this % stops a space
%\thanks{\textcopyright~2018 IEEE}% <-this % stops a space
%\thanks{}% <-this % stops a space
\thanks{\textcopyright~2018 IEEE}}

%\markboth{Control Systems Letters,~Vol.~..., December~2018}%
%{Shell \MakeLowercase{\textit{et al.}}: x}

\maketitle
\thispagestyle{empty}

\pagestyle{empty}

\begin{abstract}
It is a known fact that not all controllable systems can be asymptotically stabilized by a continuous static feedback. Several approaches have been developed throughout the last decades, including time-varying, dynamical and even discontinuous feedbacks. In the latter case, the sample-and-hold framework is widely used, in which the control input is held constant during sampling periods. Consequently, only practical stability can be achieved at best. Existing approaches often require solving optimization problems for finding stabilizing control actions exactly. In practice, each optimization routine has a finite accuracy which might influence the state convergence. This work shows, what bounds on optimization accuracy are required to achieve prescribed stability margins. Simulation studies support the claim that optimization accuracy has high influence on the state convergence.
\end{abstract}


% \fi
\IEEEpeerreviewmaketitle

\section{Introduction}\label{sec:intro}

\IEEEPARstart{I}{t} is a classical result due to \cite{Brockett1983-stabilization} which demonstrates that not every dynamical system can be stabilized by a continuous static (dependent on the state only) feedback. Since then, comprehensive work has been done in search for alternative solutions to this issue, including discontinuous feedbacks  \cite{Sontag1990-stabilization-survey,Clarke1997-stabilization,Fontes2001-MPC-SH,Fontes2003-opt-ctrl-discont,Cortes2008-discont-dyn-sys,Clarke2009-slid-mode-stab,Braun2017-SH-stabilization-Dini-aim}. This work falls specifically into the framework of sample-and-hold stabilization where the control actions are held constant during specially chosen sample time periods. 
%Such a methodology is justified from the practical standpoint by the fact that most of the controllers nowadays are implemented in a digital form. 
In this case only practical stability can be achieved at best. 
%That is, the state can be driven to an arbitrarily small vicinity of the origin. 
A concrete implementation of this setup used inf-convolutions of the given (generally non-smooth) CLF to compute its specific proximal subgradients \cite{Clarke1997-stabilization}. Using these, the control actions were calculated by computing optimizers for the CLF decay condition. In most of these techniques, in general, one has to deal with non-linear optimization problems. Unfortunately, each optimization routine has a finite accuracy which leads to approximate optimizers at best. Consequently, stability properties might be compromised. Whereas robustness properties with regards to system, input and measurement disturbance of CLF are thoroughly studied \cite{Clarke1997-stabilization,Sontag1999-stabilization-disturb,Kellett2004-Dini-aim,Kellett2000-Dini-aim,Ledyaev1997-stabilization-meas-err,Clarke2011-discont-stabilization}, additional attention should be paid to the uncertainty related to optimization inaccuracy.
%Specifically, in the inf-convolution-based feedback design, certain qualitative obstructions related to the loss of the property of being a proximal subgradient happen (see details  in the next sections).
Specifically, in inf-convolution-based feedback designs (see \eg \cite{Clarke1997-stabilization}), stabilizing control actions are computed using certain proximal subgradients. Under non-exact optimization, the property of being a proximal subgradient may be lost (details in Section \ref{sec:contributions}).
%In another technique, called Dini aiming, recently applied to the non-holonomic integrator \cite{Braun2017-SH-stabilization-Dini-aim}, a lower bound for the CLF decay is sought by optimization-based feedback.
The new result of this work is a theorem on practical stabilization by an inf-convolution-based feedback under approximate optimizers. Simulation studies of Section \ref{sec:simulation} show significant influence of optimization accuracy on the system performance. Section \ref{sec:prelim} presents some basics of the sample-and-hold framework. Section \ref{sec:contributions} summarizes the contributions, whereas Section \ref{sec:thr-res} presents the theoretical results. Discussion, importance of the made assumptions and relation to various types of uncertainty can be found in Section \ref{sec:discussion}.

\textit{Notation}: $\co(D)$ denotes the closed convex hull of a set $D$, $\ball_R(x)$ denotes a closed ball in $\R^n$ centered at $x$. If the center is the origin, the notation is simplified to just $\ball_R$. The scalar product is denoted by $\langle \bullet, \bullet \rangle$.

\section{Preliminaries}\label{sec:prelim}

Consider the following dynamical system:
\begin{equation}
	\dot x = f(x,u).
	\tag{Sys}
	\label{eqn:sys}
\end{equation}
Here, $x \in \R^n$ and $u \in \U \subseteq \R^m$ denote the system state and input, respectively. The vector field $f: \R^n \times \U \ra \R^n$ is assumed to satisfy the following local Lipschitz condition:
\begin{align*}
	& \forall c \in \R^n, r > 0, \forall x, y \in \{ z : \|z - c\| \le r \}, \forall u \in \U \\
	& \| f(x, u) - f(y, u) \| \le L_f(c, r) \cdot \| x - y \|.
	\tag{Lip}
	\label{eqn:Lip-cond} 
\end{align*} 
In general, one seeks a static feedback of the form $u = \kappa(x)$. In the sample-and-hold setup, the control is held constant during sampling periods of length $\delta$ as follows:
\begin{equation}
	\begin{array}{ll}
		\dot x 	= f(x, u_k), t \in [k \delta, (k+1) \delta], u_k \equiv \kappa(x(k \delta)).
	\end{array}
	\tag{SH}
	\label{eqn:sys-SH}
\end{equation}
It should be noted that the results in this work can be straightforwardly generalized to non-periodic partitions of the time axis with a specification of the maximal step size of $\delta$, but, for simplicity, sampling is assumed periodic here.
%In this case, the right-hand side of the system model is measurable in $t$ and the system trajectory exists in the sense of Caratheodory (for general discussion of discontinuous dynamical systems, refer \eg to \cite{Cortes2008-discont-dyn-sys,Filippov2013-discont-dyn-sys}).
The control goal is to \emph{practically stabilize} \eqref{eqn:sys-SH} in the sense that, for any given $R > r > 0$, there is a sufficiently small $\delta > 0$ such that any trajectory starting in the ball $\ball_R$ is bounded and eventually enters the ball $\ball_r$ within a time $T$ that depends uniformly on $R, r$ and the trajectory stays there after $T$.

Stabilization is achieved by utilizing a locally Lipschitz-continuous, proper and positive-definite CLF $V: \R^n \ra \R$ for which there exists a continuous function $w: \R^n \ra \R, x \ne 0 \implies w(x) > 0$ satisfying the \emph{decay condition}: for any compact set $\X \subseteq \R^n$, there exists a compact set $\U_{\X} \subseteq \U$ such that
\begin{equation}
	\forall x \in \X \inf_{\vartheta \in \co(f(x, \U_{\X}))} \; D_{\vartheta} V(x) \le -w(x).
	\tag{Dec}
	\label{eqn:decay-cond}
\end{equation} 
Here, $D_{\vartheta} V(x)$ denotes the \emph{lower directional Dini derivative}, defined as follows \cite{Clarke1997-stabilization}:
\begin{align}
	\begin{split}
	D_{\vartheta} V(x) & \triangleq \liminf_{\mu \ra 0^+} \tfrac{ V(x + \mu \vartheta) - V(x)}{\mu}. 
	\end{split}
	\tag{Dini}
	\label{eqn:Dini-der}
\end{align} 
Existence of a CLF with the above property is guaranteed for \emph{asymptotically controllable} systems \cite{Sontag1995-nonsmooth-CLF}. To actually practically stabilize the system in the sense of sample-and-hold as per \eqref{eqn:sys-SH}, one needs to compute the control actions $u_k, k \in \N$ from the given CLF $V$. Several techniques exist for this matter. For example, Clarke et al. \cite{Clarke1997-stabilization} used \emph{inf-convolutions} of $V$, which, for a given $0 < \alpha < 1$, are defined as follows:
\begin{equation}
	V_{\alpha}(x) := \inf_{y \in \R^n} \left( V(y) + \tfrac{1}{2\alpha^2} \| y-x \|^2 \right).
	\label{eqn:inf-conv}	
\end{equation}
Clarke et al. \cite{Clarke1997-stabilization} then used a minimizer $y_{\alpha}(x)$ of \eqref{eqn:inf-conv} for each given $x$ to determine the vector
\begin{equation}
	\zeta_{\alpha}(x) := \tfrac{x - y_{\alpha}(x)}{\alpha^2} 
	\label{eqn:subgrad-alpha}
\end{equation}
which also happens to be a particular \emph{proximal subgradient} of $V$ at the respective $y_{\alpha}(x)$ as follows:   
\begin{equation}
	V(z) \ge V(y_{\alpha}(x)) + \langle \zeta_{\alpha}(x) , z-y_{\alpha}(x) \rangle - \tfrac{1}{2 \alpha^2} \|z-y_{\alpha}(x)\|^2 
	\label{eqn:subgrad-alpha-cond}
\end{equation}
for every $z \in \R^n$ (for a comprehensive study on non-smooth analysis and various notions involved in this work, the reader should refer \eg to \cite{Clarke2008-nonsmooth-analys}). Suppose that $y_{\alpha}(x)$ lies within some compact set $\Y$. Then, the control actions are determined by finding minimizers for
\begin{equation}
	\inf_{u \in \U_{\Y}} \langle \zeta_{\alpha}(x), f(y_{\alpha}(x),u) \rangle
	\label{eqn:ctrl-action-OP}
\end{equation}
so as to satisfy the decay condition \eqref{eqn:decay-cond} using the fact that, for any proximal subgradient $\zeta$ of $V$ at $x$ and any vector $\vartheta$, it holds that
\begin{equation}
	\langle \zeta,\vartheta \rangle \le D_{\vartheta} V(x).
	\label{eqn:subdiff-and-Dini}
\end{equation}
Another approach is based on a technique called \emph{Dini aiming} \cite{Kellett2000-Dini-aim,Kellett2004-Dini-aim}, with a recent application to the non-holonomic integrator \cite{Braun2017-SH-stabilization-Dini-aim}, by solving special optimization problems involving the CLF and the system dynamics model $f$. Both the inf-convolution and Dini aiming methods require solving optimization problems with the respective analyzes done under the assumption that the optimization problems be solved exactly (for details, please refer \eg to \cite{Clarke1997-stabilization,Braun2017-SH-stabilization-Dini-aim}). 

In the current work, it is investigated if practical stabilization can be achieved if only approximate minimizers for \eqref{eqn:inf-conv} and \eqref{eqn:ctrl-action-OP} can be found. The scope of the current work is focused on the inf-convolution technique as the basis, whereas the developed methodology could be used for Dini aiming in a future study. The summary is given in the next section.

\section{Work aim and contributions}\label{sec:contributions}

The stabilizing feedback by inf-convolution amounts to solving the optimization problem \eqref{eqn:inf-conv} first, then \eqref{eqn:ctrl-action-OP} at every sample state $x$. Such a feedback will be called \emph{InfC-feedback} from now on. The main concern of this work is to investigate robustness properties of the CLF under an approximate InfC-feedback in the following sense: for a fixed $x$, \eqref{eqn:inf-conv} and \eqref{eqn:ctrl-action-OP} are only solved approximately. That is, an approximate $y_{\alpha}^{\eps}(x)$ is found so that:
\begin{equation}
	V(y_{\alpha}^{\eps}(x)) + \tfrac{1}{2\alpha^2} \| y_{\alpha}^{\eps}(x)-x \|^2 \le V_{\alpha}(x) + \eps_x,
\label{eqn:y-eps-def}
\end{equation}
where $\eps_x > 0$ characterizes the accuracy of the numerical optimization routine and, in general, may depend on $x$. Consequently, an approximate control action $\kappa^{\eta}_x$ is found so as to satisfy:
\begin{equation}
	\langle \zeta_{\alpha}^{\eps}(x), f(y_{\alpha}^{\eps}(x), \kappa^{\eta}_x) \rangle \le \inf_{u \in \U_{\Y}} \langle \zeta_{\alpha}^{\eps}(x), f(y_{\alpha}^{\eps}(x), u) \rangle + \eta_x,
	\label{eqn:ctrl-action-eps-OP}
\end{equation}  
where
\begin{equation}
	\zeta_{\alpha}^{\eps}(x) := \tfrac{x-y_{\alpha}^{\eps}(x)}{\alpha^2} 
	\label{eqn:subgrad-alpha-eps}
\end{equation}
is a \emph{proximal} $\eps_x$-\emph{subgradient}, and $\eta_x$, like $\eps_x$, is related to the optimization accuracy. The following notion summarizes practical stabilization under approximate optimizers:
\begin{definition}
	An InfC-feedback is said to practically stabilize \eqref{eqn:sys} in the sense of sample-and-hold \eqref{eqn:sys-SH} under approximate optimizers if, for any data $R > r > 0$, there is a sufficiently small $\delta > 0$ such that any closed-loop trajectory $x(t), t \ge 0, x(0) \in \ball_R$ is bounded and enters and stays in the ball $\ball_r$ within a time $T$ depending uniformly on the data $R, r$, provided that, at every sample state $x_k:= x(\delta k), k \in \N$, the accuracies $\eps_{x_k}$ and $\eta_{x_k}$ in \eqref{eqn:y-eps-def} and \eqref{eqn:ctrl-action-eps-OP}, accordingly, are sufficiently small.
	\label{def:pract-stab-approx-opt} 
\end{definition} 		
The \textbf{major problem} is that, in general, $\zeta_{\alpha}^{\eps}(x)$ is \textbf{not even a proximal subgradient} of $V$ at $y_{\alpha}^{\eps}(x)$. As a result, the relation \eqref{eqn:subdiff-and-Dini} does not apply for $\zeta_{\alpha}^{\eps}(x)$ in place of $\zeta$.
It might have been attractive to use some continuity properties of proximal subdifferentials \ie the sets of all proximal subgradients. Continuity properties of ordinary (not proximal) subdifferentials were studied by \cite{Gregory1980-subdiff-cont,Zalinescu2007-subdiff-cont}. However, these results cannot be directly applied to proximal subdifferentials. A proximal subdifferential may happen to be empty, but there still may exist proximal $\eps$-subgradients. To overcome the described difficulties, the following is assumed in the current work:
\begin{assumption}
	The CLF $V$ is globally lower Dini-differentiable and the $\liminf$ in \eqref{eqn:Dini-der} is locally uniform. That is, for any compact sets $\Y, \F \in \R^n$, and any $\nu > 0$, there exists $\mu > 0$ such that, for any $y \in \Y, \vartheta \in \F, 0 < \mu' \le \mu$, it holds that $\Big| \tfrac{ V(y + \mu' \vartheta) - V(y)}{\mu'} - D_{\vartheta}V(y) \Big| \le \nu$.
	\label{asm:uniform-diffty}
\end{assumption}

The main \textbf{contribution} of this work is to investigate practical stabilization under approximate optimizers in the sense of Definition \ref{def:pract-stab-approx-opt} under Assumption \ref{asm:uniform-diffty}. The main result is summarized in the next section. Assumption \ref{asm:uniform-diffty} as well as further robustness considerations are discussed in Section \ref{sec:discussion}. 

\section{Theoretical results}\label{sec:thr-res}

The following theorem summarizes the major result on practical stabilization by InfC-feedbacks under approximate optimizers.

\begin{theorem}
	Consider the control system in the sample-and-hold format \eqref{eqn:sys-SH}. Assume that $f$ satisfies the Lipschitz condition \eqref{eqn:Lip-cond} and there exists a CLF $V$ with a decay condition \eqref{eqn:decay-cond} satisfying Assumption \ref{asm:uniform-diffty}. Then, the system \eqref{eqn:sys-SH} can be practically stabilized by an InfC-feedback in the sense of \eqref{eqn:sys-SH} under approximate optimizers.
	\label{thm:pract-stab}
\end{theorem}
\begin{IEEEproof}
	The proof is divided into three steps. The first step deals with necessary preparations and parameter determination. In the second step, a relaxed decay condition of the inf-convolution of the CLF is shown. The actual decay is demonstrated in the third step along with the estimation of the reaching time $T$. The proof uses technical Lemmas \ref{lem:y-eps-near}, \ref{lem:V-between} which are found in the appendix. 
	
	\textbf{Step 1. Preliminaries}
	
	Let $0 < r < R$ be the radii of the target and starting ball, respectively. First, 	two non-decreasing functions, $\rho_{V}$ and $\lambda_{V}$, with the properties
	\begin{align*}
		\forall x \in \R^n, r, v > 0 \spc \; \spc & V(x) \le \rho_{V}(r) \implies \|x\| \le r, \\
		& V(x) \ge v \implies \lambda_V(v) \le \|x\|
	\end{align*}
	are constructed. Due to Lemma 3.5 in \cite{Khalil1996-nonlin-sys}, there exist two class $\mathcal K$ functions $\alpha_1, \alpha_2$ with the properties:
	\begin{align*}
		\forall x \in R^n \spc \alpha_1(\|x\|) \le V(x) \le \alpha_2(\|x\|).
	\end{align*}
	Then, simply taking $\rho_V(r) := \alpha_1(r)$, it holds that
	\begin{align*}
		\alpha_1(\|x\|) \le V(x) \le \alpha_1(r) \implies \|x\| \le r. 
	\end{align*}
	For $\lambda_V(r)$, take $\alpha_2^{-1}$ which is non-decreasing as well. Then,
	\begin{align*}
		\alpha_2(\|x\|) \ge V(x) \ge v \implies \|x\| \ge \alpha_2^{-1}(v).
	\end{align*}
	In the forthcoming step, the following bounds of $V$ in terms of its inf-convolutions using Lemma \ref{lem:V-between} from Appendix will be used:
	\begin{equation}
		V_{\alpha}(x) \le V(x) \le V_{\alpha}(x) + \eps_1.
		\label{eqn:V-between}
	\end{equation}
	Here, $\eps_1$ is a positive number that will be determined later. With this in mind, extend the use of $\rho_V(r)$ and $\lambda_V(r)$ to $V_{\alpha}$ as follows:
	\begin{align*}
	\begin{aligned}
		& V_{\alpha}(x) \le \rho_{V}(r) - \eps_1 \implies V(x) \le \rho_{V}(r) \implies \|x\| \le r, \\
		& V_{\alpha}(x) \ge v \implies V(x) \ge v \implies \|x\| \ge \lambda_{V}(v).
	\end{aligned}
	\end{align*}
	Set $\bar V := \sup_{\|x\| \le R} V(x)$. Choose $R^* > R$ with the corresponding $V^* := \sup_{\|x\| \le R^*} V(x)$ such that $\Theta := \rho_V(R^*) > \bar V$. Further, set $v^* := \rho_V(r)$ and $r^* := \lambda_V\left( \tfrac{v^*}{4} \right)$. Observe that $r^* \le r$. To see this, assume, on contrary, that $r^* > r$. Since $V$ is continuous, there exists $x$ such that $\tfrac{v^*}{4} \le V(x) \le v^*$. By the virtue of $\rho_V$, it follows that $\|x\| \le r$, whereas $\|x\| \ge \lambda_V \left( \tfrac{v^*}{4} \right) > r$, a contradiction. Now, let $\U^*$ be the compact set corresponding to the ball $\ball_{R^* + \sqrt{2 V^*}}$ in the decay condition \eqref{eqn:decay-cond}:
	\begin{align*}
		\forall x \in \ball_{R^* + \sqrt{2 V^*}} \inf_{\vartheta \in \co(f(x, \U^*))} D_{\vartheta} V(x) \le -w(x).
	\end{align*}
	Let $L_f$ be the Lipschitz constant of $f$ as per \eqref{eqn:Lip-cond} on the ball $\ball_{R^* + \sqrt{2 V^*}}$.
	Set
	\begin{align*}
		\bar f := \sup_{\substack{x \in \ball_{R^* + \sqrt{2 \bar V^*}} \\ u \in \U^*}} f(x,u)
	\end{align*}
	and determine the minimal decay rate as follows:
	\begin{align*}
		\bar w := \inf_{\tfrac{r^*}{2} \le \|x\| \le R^* + \sqrt{2 V^*}} w(x).
	\end{align*}
	Finally, let $\omega_V$ be the continuity modulus of $V$ on $\ball_{R^* + \sqrt{2 V^*}}$. That is, for all $\eps_1 > 0$, and $x, y \in \ball_{R^* + \sqrt{2 V^*}}$, it holds that
	\begin{align*}
		& \| x - y \| \le \omega_V(\eps_1) \implies |V(x) - V(y)| \le \eps_1.
	\end{align*}	
	
	\textbf{Step 2. Relaxed decay condition}
	
	A relaxed decay condition is established by addressing the scalar product $\langle \zeta_{\alpha}^{\eps^2}(x), f(x, \kappa^{\eta}_x) \rangle$, where $\zeta_{\alpha}^{\eps^2}(x)$ will be a proximal $\eps_x^2$-subgradient with $\eps_x$ determined later on sample-wise at $x_k, k \in N$, and $\kappa^{\eta}_x$ will be the control action determined as an approximate optimizer of $\inf_{u \in \U^*} \langle \zeta_{\alpha}^{\eps^2}(x), f(x, u) \rangle$ with the yet-to-be-determined accuracy $\eta_x$. 
	
	Using the Lipschitz constant $L_f$ of $f$, it holds, for all $x$ with $r^* \le \|x\| \le R^*$, that: 
	\begin{equation*}
		\begin{split}
			& \langle \zeta_{\alpha}^{\eps^2}(x), f(x, \kappa^{\eta}_x) \rangle \le \eta_{x} \\
			& +\underbrace{\inf_{u \in \U^*} \langle \zeta_{\alpha}^{\eps^2}(x), f(y_{\alpha}^{\eps^2}(x), u) \rangle}_{\circled{1}} + L_f \underbrace{\| \zeta_{\alpha}^{\eps^2}(x) \| \cdot \|y_{\alpha}^{\eps^2}(x)-x\|}_{\circled{2}},
		\end{split}	
	\end{equation*}	
	where $y_{\alpha}^{\eps^2}(x)$ is an $\eps_x^2$-minimizer for the respective inf-convolution \eqref{eqn:inf-conv}. Notice $\eta_x$ in the right-hand side which indicates that the control actions $\kappa^{\eta}_x$ are found only in an approximate format.
	
	Now, examine the term \circled{1}. 
	This is the term which the approximate proximal subgradient condition in Step 3 will be applied to.
	But first, it requires a fixed $\alpha$.  Therefore, $\eps$ will be determined later when all the conditions on $\alpha$ are established. Regarding the term \circled{2}, observe that, by definitions of $y_{\alpha}^{\eps^2}(x)$ and $\zeta_{\alpha}^{\eps^2}(x)$ according to \eqref{eqn:y-eps-def}, \eqref{eqn:subgrad-alpha-eps}, respectively, and Lemma \ref{lem:y-eps-near} from Appendix,
	\begin{align*}
		\| \zeta_{\alpha}^{\eps^2}(x) \| \cdot \| y_{\alpha}^{\eps^2}(x) - x \| = & \tfrac{1}{\alpha^2} \| y_{\alpha}^{\eps^2}(x) - x \|^{2} \\
		\le & 2(V(x) - V(y_{\alpha}^{\eps^2})) + 2 \eps_x^2.
	\end{align*}	
	If $\alpha$ is chosen to satisfy the following condition:
	\begin{align}
		& \sqrt{2 V^*} \alpha \le \min \left\{ 1, \tfrac{r^*}{2}, \omega_V(\eps_1) \right\},
		\label{eqn:alpha-1st-cond}
	\end{align}	
	where $\eps_1$ is a positive number yet to be determined, then, it holds that
	\begin{equation}
		\begin{aligned}
			 \langle \zeta_{\alpha}^{\eps^2}(x), f(x, \kappa^{\eta}_x) \rangle \le & \inf_{u \in \U^*} \langle 	\zeta_{\alpha}^{\eps^2}(x), f(y_{\alpha}^{\eps^2}(x), u) \rangle \\
			 & + 2 L_f ( \eps_1 + \eps_x^2 ) + \eta_x.
		\end{aligned}		 
		\label{eqn:relaxed-decay-cond-interm}
	\end{equation}
	
	Also, observe that $\tfrac{r^*}{2} \le \| y_{\alpha}^{\eps^2}(x) \| \le R^* + \sqrt{2 V^*}$ and \eqref{eqn:V-between} holds. The parameters $\eps_x, \eps_1, \eta_x$ in \eqref{eqn:relaxed-decay-cond-interm} will be determined in Step 3.
	
	\textbf{Step 3. Decay}
	
	First, it follows directly from the definitions, that a variant of ``Taylor expansion'' for $V_{\alpha}$ holds in an approximate format \ie for all $h, \vartheta \in \R^n$,
	\begin{align}
		 V_{\alpha}(x + h \vartheta) \le V_{\alpha}(x) + h \langle \zeta_{\alpha}^{\eps^2}(x) , \vartheta \rangle + \tfrac{h^2 \|\vartheta\|^2}{2\alpha^2} + \eps_x^2.
		\label{eqn:Taylor-expan}
	\end{align}
	Proceed by induction over the sample time periods. Suppose that the trajectory of \eqref{eqn:sys-SH} exists locally at the time period $[k \delta, (k+1) \delta]$. Assume $V_{\alpha}(x_k) \le \bar V$. To show that the trajectory exists on the entire sample time period, observe that, due to \eqref{eqn:V-between}, for $t \in [k \delta, (k+1) \delta]$,
	\begin{equation}
		V(x(t)) \le V_{\alpha}(x(t)) + \eps_1 \le \bar V + \eps_1.
		\label{eqn:eps1-1st-cond}
	\end{equation}
	Therefore, $\eps_1$ can be chosen so as to satisfy the following condition $V(x(t)) \le \Theta$. Consequently, $x(t) \in \ball_{R^*}$ and thus the overshoot is bounded. In the following, it is shown that $V_{\alpha}(x_k)$ may only decay down to a prescribed limit sample-wise \ie for $k \in \N \spc V_{\alpha}(x_{k+1}) < V_{\alpha}(x_k)$ until $V_{\alpha}(x_k) \le v^*$. Therefore, boundedness of the trajectory on each sample period is secured. The following cases are now possible.
	
	\emph{Case 1}: $V_{\alpha}(x_k) \ge \tfrac{v^*}{2}$. Use the ``Taylor expansion'' \eqref{eqn:Taylor-expan} to deduce, for any $t \in [k \delta, (k+1)\delta], \Delta t := t- k \delta$, that
	\begin{equation}
		\begin{split}
		V_{\alpha}(x(t)) - V_{\alpha}(x_k) = & V_{\alpha}(x_k + \delta F_k) - V_{\alpha}(x_k) \le \\
		& \delta \langle \zeta_{\alpha}^{\eps^2}(x_k) , F_k \rangle + \tfrac{\delta^2 \|F_k\|^2}{2 \alpha^2} + \eps^2_{x_k},
		\end{split}
		\label{eqn:CLF-bn-nodes}
	\end{equation}
	where $F_k$ is defined by the integral form of the trajectory:
	\begin{align*}
		x(t) = & x_k + \int \limits_{k \delta}^{t} f(x(\tau), \kappa^{\eta}_{x_k} ) \, \diff\tau = \\
		& x_k + \delta \underbrace{\left( \tfrac{1}{\delta} \int \limits_{k \delta}^{t} f(x(\tau), \kappa^{\eta}_{x_k}) \, \diff\tau \right) }_{=: F_k}
	\end{align*}
	and, with the property $\|F_k\| \le \tfrac{1}{\delta} \Delta t \bar f$, $F_k$ can be re-expressed as 
	\begin{align*}
		& F_k = \tfrac{\Delta t}{\delta} f(x_k,\kappa^{\eta}_{x_k}) + \underbrace{ \tfrac{1}{\delta} \int \limits_{k \delta}^{t} \left( f(x(\tau),\kappa^{\eta}_{x_k}) - f(x_k,\kappa^{\eta}_{x_k}) \right) \, \diff\tau }_{=:A}.	
	\end{align*}
	Now, since $\|A\| \le \tfrac{\Delta t^2}{\delta}L_f \bar f$, the scalar product $\langle \zeta_{\alpha}^{\eps^2}(x_k) , F_k \rangle$, is bounded as follows:
	\begin{align*}
		& \langle \zeta_{\alpha}^{\eps^2}(x_k) , F_k \rangle \le \\ 
		& \tfrac{\Delta t}{\delta} \langle \zeta_{\alpha}^{\eps^2}(x_k) , f(x_k,\kappa^{\eta}_{x_k}) \rangle + \| \zeta_{\alpha}^{\eps^2}(x_k) \| \tfrac{\Delta t^2}{\delta}L_f \bar f \le \\
		& \tfrac{\Delta t}{\delta} \left( \inf_{u \in \U^*} \langle \zeta_{\alpha}^{\eps^2}(x_k), f(y_{\alpha}^{\eps^2}(x_k), u) \rangle + 2 L_f ( \eps_1 + \eps^2_{x_k} ) + \eta_{x_k} \right) \\
		& + \tfrac{\sqrt{2 V^*}}{\alpha} \tfrac{\Delta t^2}{\delta}L_f \bar f,
	\end{align*}	
	which follows from the relaxed decay condition \eqref{eqn:relaxed-decay-cond-interm}. Therefore,
	\begin{equation}
		\begin{aligned}
			& V_{\alpha}(x(t)) - V_{\alpha}(x_k) \le \Delta t \Big( \inf_{u \in \U^*} \langle \zeta_{\alpha}^{\eps^2}(x_k), f(y_{\alpha}^{\eps^2}(x_k), u) \rangle \\ 
			& + 2 L_f ( \eps_1 + \eps_{x_k}^2 ) + \eta_{x_k} + \tfrac{\sqrt{2 V^*}}{\alpha} \Delta t L_f \bar f + \tfrac{\Delta t \bar f^2}{2 \alpha^2} \Big) + \eps_{x_k}^2.
		\end{aligned}		
	\end{equation}
	In particular, for $t=(k+1) \delta$, observe that
	\begin{equation}
		\begin{aligned}
			& V_{\alpha}(x_{k+1}) - V_{\alpha}(x_k) \le \delta \Big( \inf_{u \in \U^*} \langle \zeta_{\alpha}^{\eps^2}(x_k), f(y_{\alpha}^{\eps^2}(x_k), u) \rangle \\ 
			& + 2 L_f ( \eps_1 + \eps^2_{x_k} ) + \eta_{x_k} + \tfrac{\sqrt{2 V^*}}{\alpha} \delta L_f \bar f + \tfrac{\delta \bar f^2}{2 \alpha^2} \Big) + \eps^2_{x_k}.
		\end{aligned}	
		\label{eqn:relaxed-decay-cond-2nd-interm}	
	\end{equation}	
	\emph{Case 2}: $V_{\alpha}(x_k) \le \tfrac{3 v^*}{4}$. If the sample period size $\delta$ satisfies $\delta \bar f \le \omega_V (\eps_2)$ for some $\eps_2 > 0$, then
	\begin{align}
		V_{\alpha}(x(t)) \le V_{\alpha}(x_k) + \eps_2. 
	\end{align}	 
	Choosing $\eps_2 \le \tfrac{v^*}{8}$ guarantees that $V_{\alpha}(x(t)) \le \tfrac{7 v^*}{8}$. Now, recall the condition \eqref{eqn:V-between} where $\eps_1$ needed to be determined. Taking into account \eqref{eqn:eps1-1st-cond}, $\eps_1$ can be determined so as to satisfy the condition $V(x(t)) \le v^*$ which also implies $\| x(t) \| \le r$.
	
	Finally, the necessary parameters can be determined to establish the decay. So far, certain bounds on $\eps_1, \eps_2$ have already been established. Recalling \eqref{eqn:relaxed-decay-cond-2nd-interm}, constrain $\delta < 1$ and $\eps_1$ also via $2 L_f \eps_1 \le \tfrac{\bar w}{20}$. Notice that this indirectly constrains $\alpha$ and $\eps_{x_k}$ as well by Lemma \ref{lem:V-between} from Appendix. From now on, $\alpha$ is fixed. Bound the sample step size $\delta$ as follows:
	\begin{align}
		\eta_{x_k} \le \tfrac{\bar w}{20}, & & \tfrac{\delta \bar f^2}{2 \alpha^2} \le \tfrac{\bar w}{20}, & & \delta \tfrac{\bar w}{2} \le \tfrac{v^*}{4}.
		\label{eqn:conditions1}
	\end{align}
	Force $\delta$ to additionally satisfy $\tfrac{\sqrt{2 V^*}}{\alpha}L_f \bar f \delta \le \tfrac{\bar w}{20}$.
	From now on, $\delta$ is considered fixed. Constrain $\eps$ further by the conditions:
	\begin{align}
		& \eps^2_{x_k} \le \delta \tfrac{\bar w}{20}, & 2 L_f \eps^2_{x_k} \le \tfrac{\bar w}{20}.
	\end{align}
	Now, the most crucial part, optimization accuracy for the subgradients is addressed. In the following, it is shown how to find a further bound on $\eps^2_{x_k}$ such that the next relation holds:
	\begin{equation}
		\inf_{u \in \U^*} \langle \zeta_{\alpha}^{\eps^2}(x_k) , f(y_{\alpha}^{\eps^2}(x_k), u) \rangle \le -\tfrac{3 \bar w}{4}. 
		\label{eqn:conditions2}  
	\end{equation}
	
	\noindent \textbf{\emph{Condition on approximate proximal subgradients}}

	%%%%%%%%%%%%%%%%%%%%%%%%%%%%%%%%%%%%%%%%%%%%%%%%%%%%%%%%%%%%%%%%%%%%%%%%%%%%%%%%%%%%%%%%%%%%%%%	
	%%%%%%%%%%%%%%%%%%%%%%%%%%%%%%%%%%%%%%%%%%%%%%%%%%%%%%%%%%%%%%%%%%%%%%%%%%%%%%%%%%%%%%%%%%%%%%%
	%%%%%%%%%%%%%%%%%%%%%%%%%%%%%%%%%%%%%%%%%%%%%%%%%%%%%%%%%%%%%%%%%%%%%%%%%%%%%%%%%%%%%%%%%%%%%%%
	%%%%%%%%%%%%%%%%%%%%%%%%%%%%%%%%%%%%%%%%%%%%%%%%%%%%%%%%%%%%%%%%%%%%%%%%%%%%%%%%%%%%%%%%%%%%%%%
	%%%%%%%%%%%%%%%%%%%%%%%%%%%%%%%%%%%%%%%%%%%%%%%%%%%%%%%%%%%%%%%%%%%%%%%%%%%%%%%%%%%%%%%%%%%%%%%

	First, observe that, for any $z \in \R^n$, the following holds:
	\begin{align*}
		V(z) \ge & V(y_{\alpha}^{\eps^2}(x_k)) + \langle \zeta_{\alpha}^{\eps^2}(x_k) , z-y_{\alpha}^{\eps^2}(x_k) \rangle \\
		& - \tfrac{1}{2 \alpha^2} \|z-y_{\alpha}^{\eps^2}(x_k)\|^2 - \eps^2_{x_k}.
	\end{align*}	
	In particular, for any $\vartheta \in R^n$, the following bound applies:
	\begin{align*}
		V(y_{\alpha}^{\eps^2}(x_k) + \eps \vartheta) \ge & V(y_{\alpha}^{\eps^2}(x_k)) + \eps_{x_k} \langle \zeta_{\alpha}^{\eps^2}(x_k), \vartheta \rangle \\
		& - \tfrac{1}{2 \alpha^2} \eps^2_{x_k} \|\vartheta\|^2 - \eps^2_{x_k}.
	\end{align*}		
	Then, it follows that $\langle \zeta_{\alpha}^{\eps^2}(x), \vartheta \rangle$ is bounded from above by
	\begin{equation}
		\tfrac{V(y_{\alpha}^{\eps^2}(x_k) + \eps_{x_k} \vartheta) - V(y_{\alpha}^{\eps^2}(x_k))}{\eps_{x_k}} + \tfrac{1}{2 \alpha^2} \eps_{x_k} \|\vartheta\|^2 + \eps_{x_k}.
		\label{eqn:lemma-prf-scalar-product}
	\end{equation}
	By Assumption \ref{asm:uniform-diffty}, for all $\displaystyle y \in \Y_k^* := \{y \in \R^n: \|y - x_k\| \le (2 \bar V)^{\nicefrac{1}{2}} \alpha \}$,  $\vartheta \in \co ( f(\Y_k^*, \U^*) )$, there is a $\mu > 0$ such that
	\begin{equation}
		\Big| \tfrac{V(y) + \mu' \vartheta) - V(y)}{\mu'} - D_{\vartheta}V(y) \Big| \le \tfrac{\bar w}{5}
	\end{equation}
	holds for any $0 < \mu' \le \mu$. Notice that any approximate optimizer $y_{\alpha}^{\eps^2}(x_k)$ lies in $\Y_k^*$ disregarding $\eps^2_{x_k}$ by Lemma \ref{lem:y-eps-near} (see Appendix). Set $\eps^2_{x_k}$ not greater than $\mu$ and conclude that
	\begin{equation}
		\Big| \tfrac{V(y_{\alpha}^{\eps^2}(x_k)) + \eps^2_{x_k} \vartheta) - V(y_{\alpha}^{\eps^2}(x_k))}{\eps^2_{x_k}} - D_{\vartheta}V(y_{\alpha}^{\eps^2}(x_k)) \Big| \le \tfrac{\bar w}{5}
		\label{eqn:asm1-applied}
	\end{equation}	
	holds for any $\vartheta \in \co ( f(y_{\alpha}^{\eps^2}(x_k), \U^*) )$. Remember that $\mu$ depends on $x_k$ and, therefore, $\eps^2_{x_k}$ does so as well.
	Putting together the bound \eqref{eqn:lemma-prf-scalar-product} and relation \eqref{eqn:asm1-applied} yields
	\begin{align*}
		\langle \zeta_{\alpha}^{\eps^2}(x_k), f(y_{\alpha}^{\eps^2}(x_k), u) \rangle \le & D_{f(y_{\alpha}^{\eps^2}(x_k), u)}V(y_{\alpha}^{\eps^2}(x_k)) + \tfrac{\bar w}{5} \\
		& + \tfrac{\eps_{x_k} }{2 \alpha^2} \|f(y_{\alpha}^{\eps^2}(x_k), u)\|^2 + \eps_{x_k}
	\end{align*}
	for any $u \in \U^*$. Consequently, it holds that
	\begin{align*}
		& \inf_{\vartheta \in \co ( f( y_{\alpha}^{\eps^2}(x_k), \U^*)) } \langle \zeta_{\alpha}^{\eps^2}(x_k), \vartheta \rangle \le \\  & \inf_{\vartheta \in \co ( f( y_{\alpha}^{\eps^2}(x_k), \U^*)) } D_{\vartheta}V(y_{\alpha}^{\eps^2}(x_k)) + \tfrac{\bar w}{5} + \tfrac{\eps_{x_k}}{2 \alpha^2} \bar f^2 + \eps_{x_k} \le \\
		& -\tfrac{4 \bar w}{5} + \tfrac{\eps_{x_k}}{2 \alpha^2} \bar f^2 + \eps_{x_k}.
	\end{align*}
	Therefore, the following bound
	\begin{equation}
		\eps_{x_k} \le \tfrac{\bar w \alpha^2}{10(\bar f^2 + 2 \alpha^2)}
		\label{eqn:conditions3}
	\end{equation}	
	ensures the condition 
	\begin{equation*}
		\inf_{\vartheta \in \co ( f( y_{\alpha}^{\eps^2}(x_k), \U^*)) } \langle \zeta_{\alpha}^{\eps^2}(x_k), \vartheta \rangle \le - \tfrac{3 \bar w}{4}.
	\end{equation*}	
	By \cite[Lemma~1]{Kellett2000-Dini-aim}, it also holds that
	\begin{equation}
		\inf_{u \in \U^*} \langle \zeta_{\alpha}^{\eps^2}(x_k) , f(y_{\alpha}^{\eps^2}(x_k), u) \rangle \le -\tfrac{3 \bar w}{4},
		\label{eqn:decay-approx-subgrad}
	\end{equation}
	which is the required condition.	
	
	\noindent \textbf{\emph{Reaching time}}
	
	Putting all the relevant constraints together yields an inter-sample decay rate of $\delta \tfrac{\bar w}{2}$ on $V_{\alpha}$. The time $T$ of reaching \emph{Case 2} can be determined given by $T_{\alpha} = 4 \tfrac{\bar V - \tfrac{v^*}{2}}{\bar w}$. As soon as \emph{Case 2} is reached, two subcases in the consequent sampling periods are possible.
	
	\noindent \emph{Subcase 2.1}: $\tfrac{v^*}{2} \le V_{\alpha}(x_k) \le v^*$. \emph{Subcase 2.2}: $V_{\alpha}(x_k) \le \tfrac{3 v^*}{4}$.
	
	\noindent Being in \emph{Subcase 2.1}, $V_{\alpha}$ can either stay in it in the next sample period or jump to \emph{Subcase 2.2} since the decay condition holds. Being in \emph{Subcase 2.2}, $V_{\alpha}$ can either stay in it or move to \emph{Case 2}, from which, again, $V_{\alpha}$ can only move to one of the subcases. In both subcases, as well as in \emph{Case 2}, the trajectory stays in the target ball.	
\end{IEEEproof}
\begin{remark}
	Explicit bounds on $\delta, \alpha, \eta_{x_k}, \eps_{x_k}$ are all in principle computable from various conditions given in the proof -- in particular, via \eqref{eqn:conditions1}--\eqref{eqn:conditions2}, \eqref{eqn:asm1-applied}, \eqref{eqn:conditions3} etc.
\end{remark}
\begin{remark}
	The derived result is semi-global in the sense of \cite{Clarke1997-stabilization}. In case of a local CLF with the decay property \eqref{eqn:decay-cond} where some compact set $\Omega$ is used instead of $\R^n$, the only difference is that the starting ball $\ball_R$ must be restricted so that the corresponding $\ball_{R^* + \sqrt{2 V^*}}$ is within $\Omega$.
\end{remark}
\begin{remark}
	Stabilization can also be achieved if Assumption \ref{asm:uniform-diffty} does not hold globally \ie if $\Y$ may merely belong to some subset $\mathbb A$ of $\R^n$ such that $\mathbb A_0 := \R^n \setminus \mathbb A$ has measure zero, as long as, for each $x_k \notin \ball_{\nicefrac{r^*}{2}}$, the set $\Y^*_k$ does not lie within some set $\mathbb A_{0_{\chi}} := \{y \in \R^n: \|y - \mathbb A_0 \| \le \chi \}$, where $\chi > 0$ may be arbitrary. Otherwise, Assumption \ref{asm:uniform-diffty} can be relaxed to the following: 
	for any $\delta' > 0$, any compact sets $\Y, \F \subset \R^n$, any $\nu, \sigma >0$, there exist $\chi >0$, a set $\Y' \subseteq \Y$ and $\mu >0$ such that
	\begin{enumerate}
	\item[1)] for all $y \in \Y'$, all $\vartheta \in \F$ and all $0 < \mu' \leq \mu$, it holds that $\Big| \tfrac{ V(y + \mu' \vartheta) - V(y)}{\mu'} - D_{\vartheta}V(y) \Big| \le \nu$;
	\item[2)] there exists $M \in \N$, such that, for all $\tilde y \in \tilde{\Y} := (\Y \setminus \Y')_{\chi}$, there exists $\{\pi_i\}_{i=1}^M \in \U^M$ such that the trajectory of \eqref{eqn:sys-SH} with $\delta = \delta'$, the initial condition $\tilde y$ and under the control sequence $\{\pi_i\}_{i=1}^M$ satisfies either of the following:
		\begin{enumerate}
		\item[i)] $x(\delta' M) \in \Y \setminus \tilde \Y $ and $V(x(\delta' M)) \leq V(\tilde y) + \sigma$,
		\item[ii)] $V(x(\delta' M)) < V(\tilde{y})$.
		\end{enumerate}
	\end{enumerate}	 
	In this case, setting $\sigma \le \nicefrac{\delta \bar w}{4}, \alpha \le \nicefrac{\chi}{\sqrt{2 \bar V}}$, the local controller $\{\pi_i\}_{i=1}^M$ may be invoked whenever $\|x_k - \Y \setminus \Y'\| \le \chi$.
	\label{rem:A1-relaxation}
\end{remark}

\section{Discussion}\label{sec:discussion}

Theorem \ref{thm:pract-stab} demonstrates robustness properties of a given CLF against non-exact optimization in the InfC-feedback represented via the problems \eqref{eqn:inf-conv} and \eqref{eqn:ctrl-action-OP}. Further types of uncertainties can be introduced into the setup without compromising the practical stability provided that they are (essentially) bounded in a certain way.
For instance, an uncertainty of the type $\dot x = f(x, u) + g(t)$ can be addressed mainly by taking care of the ``Taylor expansion'' \eqref{eqn:CLF-bn-nodes}, where a respective term of the form $\int_{k \delta}^t g(\tau) \diff \tau$ appears, and bounding $g$ sufficiently along the lines of \cite{Clarke1997-stabilization}. An actuator uncertainty of the type $\dot x = f(x, u + d(t))$ can be converted into a one of the type $\dot x = f(x, u) + g(t)$ using continuity of $f$.
In this sense, an approximate InfC-feedback in the spirit of Definition \ref{def:pract-stab-approx-opt} can be made robust against the said uncertainties.
A natural question is whether the optimization accuracies $\eta_x$ and $\eps_x$ can be ``merged'' with the uncertainty $g(t)$. Whereas it is relatively simple with $\eta_x$, introducing $\eps_x$ into $g(t)$ would require special attention. The problem is that, in presence of optimization inaccuracy, as pointed out earlier, a qualitative obstruction, that is loss of the property of being a proximal subgradient and, subsequently, failing to satisfy \eqref{eqn:subdiff-and-Dini}, may occur. The parameter $\eps_x$ occurs at many places in the proof of Theorem \ref{thm:pract-stab}, but mainly in the \textbf{\emph{condition on approximate proximal subgradients}}. To derive a necessary bound on $g$ so as to address the optimization inaccuracy would require tracking back the derivations in the said fragments of the theorem. Therefore, the present result should be seen as a complementary one to the works based on InfC-feedbacks such as in \cite{Clarke1997-stabilization}.
As other types of discontinuous feedbacks, an approximate InfC-feedback is still vulnerable under a measurement error of the type $\dot x = f(x, (u(x+d(t))))$. In this regard, the presented result might be merged with the one in \cite{Ledyaev1997-stabilization-meas-err} and \cite[Theorem~E]{Sontag1999-stabilization-disturb}. The first one suggests an internal tracking controller to amend the main one, whereas the latter uses, besides an upper, a lower bound on the sampling interval to compensate for measurement noise.
Without Assumption \ref{asm:uniform-diffty}, no locally uniform bound on $\eps_x$ could be derived in Theorem \ref{thm:pract-stab}. In absence of local uniformness in the sense of Assumption \ref{asm:uniform-diffty} or its relaxation in Remark \ref{rem:A1-relaxation}, the accuracy $\eps_x$ would, in particular, depend on the direction, which is actually defined afterwards using the said accuracy, a circular reasoning. Notice an important difference to \cite{Clarke1997-stabilization} where no such an assumption was necessary. Still, Assumption \ref{asm:uniform-diffty} is merely sufficient, whereas its necessity must be further investigated. Finally, whereas $\delta$ is an absolutely crucial parameter in the analysis, the accuracy bounds $\eps_x$ and $\eta_x$ should be paid high attention as well.

\section{Simulation study}\label{sec:simulation}

\begin{figure}[h]
	\centering
	\includegraphics[width=0.8\columnwidth]{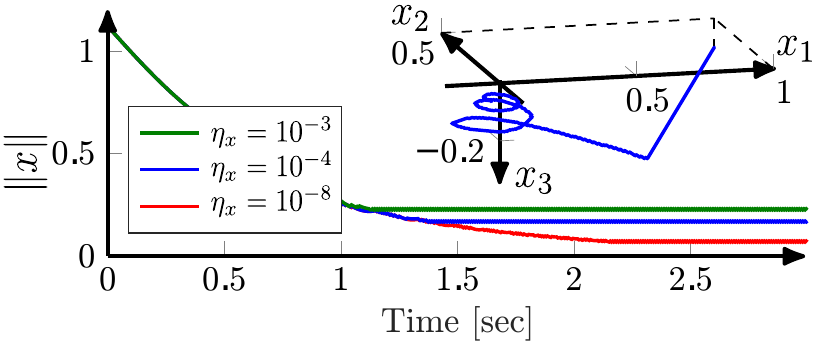}
	\caption{State norm under various $\eta_x$ and state trajectory for $\eta_x = 10^{-8}$}
	\label{fig:state-norm}
	\vspace{-6pt}
\end{figure}

\begin{figure}[h]
	\centering
	\includegraphics[width=0.8\columnwidth]{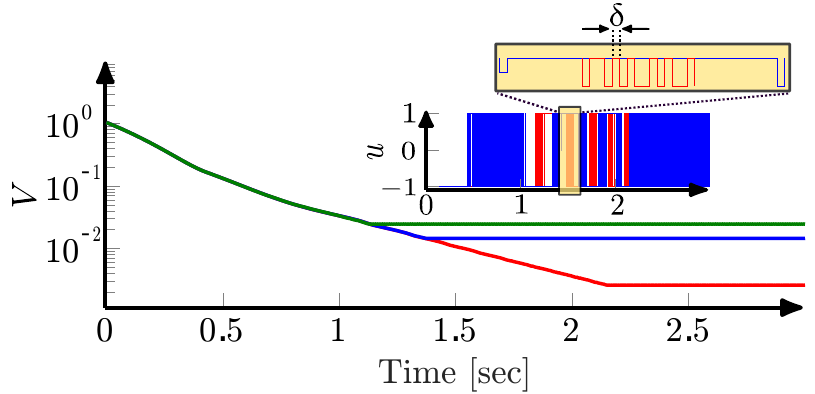}
	\caption{Behavior of $V$. The input is shown for $\eta_x = 10^{-8}$ ($u_1$\textcolor{red}{--}, $u_2$\textcolor{blue}{--} )}
	\vspace{-6pt}
	\label{fig:CLF}
\end{figure}

As the basis for the simulation study, practical stabilization of the nonholonomic integrator
\begin{align*}
(\dot{x}_1, \dot{x}_2, \dot{x}_3)^\top = (u_1, u_2, x_1 u_2 - x_2 u_1)^\top
\end{align*} 
is used. The nonholonomic integrator can be regarded as a three-wheel robot with a speed and steering control and serves thus as an important basis of a class of systems describing wheeled machines.
%Another recent work in this direction was made in \cite{Braun2017-SH-stabilization-Dini-aim}, which used Dini aiming. 
The purpose of the current study is not to compare the InfC-feedbacks to other stabilization techniques, but rather to demonstrate the effects of the optimization accuracy on stability properties.
It was shown in \cite{Braun2017-SH-stabilization-Dini-aim} that the following function is a global CLF for the nonholonomic integrator:
\begin{align*}
	V(x) = x_1^2 + x_2^2 + 2x_3^2 - 2|x_3|\sqrt{x_1^2 + x_2^2}, 
\end{align*}
under the constraint $u \in [-1 \; \, 1]^2$. It can be verified that $V$ satisfies Assumption \ref{asm:uniform-diffty} everywhere except for $\{x \in R^n: x_3=0\}_{\chi} \bigcup \{x \in R^n: x_1=0 \land x_2=0\}_{\chi}$ for any arbitrary, but fixed, $\chi > 0$. In the current simulation result, the trajectory stayed apart from $\{x \in R^n: x_3=0\} \bigcup \{x \in R^n: x_1=0 \land x_2=0\}$ (see Fig. \ref{fig:state-norm}).  
A sampling time $\delta = 0.005$ and an initial condition $x_0 = (1 , \, 0.5 , \, -0.1 )^\top$ are assumed. The parameter $\alpha$ is set equal $0.1$. The optimization problems \eqref{eqn:inf-conv}, \eqref{eqn:ctrl-action-eps-OP} of the InfC-feedback are solved until the accuracies $\eps^2_x$ and, respectively, $\eta_x$ are achieved. For simplicity, $\eps^2_x$ and $\eta_x$ are set equal. Therefore, only $\eta_x$ is specified from now on. Fig. \ref{fig:state-norm} shows the norm of the system trajectory under three different values of $\eta_x$.
It can be observed that $\eta_x$ significantly influences the stability margins, especially the vicinity into which the state trajectory converges. Fig. \ref{fig:CLF} shows the respective behaviors of $V$. A particular input signal for the case with $\eta_x = 10^{-8}$ is shown in the right-upper corner of that figure.
One can see a characteristic chattering of the control which was also observed by \cite{Braun2017-SH-stabilization-Dini-aim}. Future investigation is required for other techniques, such as Dini aiming, and also for further relaxations of Assumption \ref{asm:uniform-diffty}.

\section{Conclusion}

This work is concerned with practical stabilization of non-linear dynamical systems in the sample-and-hold framework. The key new result is an analysis of practical stability under approximate optimizers for various optimization problems. Bounds on optimization accuracy to achieve prescribed stability margins are derived. A simulation study showed significant effects of optimization accuracy on stability properties.

\bibliographystyle{plain}
\bibliography{bibl}

\begin{thebibliography}{10}

\bibitem{Braun2017-SH-stabilization-Dini-aim}
P.~Braun, L.~Gr{\"u}ne, and C.~Kellett.
\newblock Feedback design using nonsmooth control lyapunov functions: A
  numerical case study for the nonholonomic integrator.
\newblock In {\em Proceedings of the 56th IEEE Conference on Decision and
  Control}. IEEE, 2017.

\bibitem{Brockett1983-stabilization}
R.~Brockett.
\newblock Asymptotic stability and feedback stabilization.
\newblock {\em Differential geometric control theory}, 27(1):181--191, 1983.

\bibitem{Clarke2011-discont-stabilization}
F.~Clarke.
\newblock Lyapunov functions and discontinuous stabilizing feedback.
\newblock {\em {A}nnual {R}eviews in {C}ontrol}, 35(1):13--33, 2011.

\bibitem{Clarke1997-stabilization}
F.~Clarke, Y.~Ledyaev, E.~Sontag, and A.~Subbotin.
\newblock Asymptotic controllability implies feedback stabilization.
\newblock {\em IEEE Transactions on Automatic Control}, 42(10):1394--1407,
  1997.

\bibitem{Clarke2008-nonsmooth-analys}
F.~Clarke, Y.~Ledyaev, R.~Stern, and P.~Wolenski.
\newblock {\em Nonsmooth Analysis and Control Theory}, volume 178.
\newblock Springer Science \& Business Media, 2008.

\bibitem{Clarke2009-slid-mode-stab}
F.~Clarke and R.~Vinter.
\newblock Stability analysis of sliding-mode feedback control.
\newblock {\em {C}ontrol and {C}ybernetics}, 4(38):1169--1192, 2009.

\bibitem{Cortes2008-discont-dyn-sys}
J.~Cortes.
\newblock Discontinuous dynamical systems.
\newblock {\em IEEE {C}ontrol {C}ystems}, 28(3), 2008.

\bibitem{Fontes2001-MPC-SH}
F.~Fontes.
\newblock A general framework to design stabilizing nonlinear model predictive
  controllers.
\newblock {\em Systems \& Control Letters}, 42(2):127--143, 2001.

\bibitem{Fontes2003-opt-ctrl-discont}
F.~Fontes.
\newblock Discontinuous feedbacks, discontinuous optimal controls, and
  continuous-time model predictive control.
\newblock {\em International Journal of Robust and Nonlinear Control},
  13(3-4):191--209, 2003.

\bibitem{Gregory1980-subdiff-cont}
D.~Gregory.
\newblock Upper semicontinuity of subdifferential mappings.
\newblock {\em Canadian Mathematical Bulletin}, 23:11--19, 1980.

\bibitem{Kellett2004-Dini-aim}
C.~Kellett, H~Shim, and A.~Teel.
\newblock Further results on robustness of (possibly discontinuous) sample and
  hold feedback.
\newblock {\em IEEE Transactions on Automatic Control}, 49(7):1081--1089, 2004.

\bibitem{Kellett2000-Dini-aim}
C.~Kellett and A.~Teel.
\newblock Uniform asymptotic controllability to a set implies locally
  {L}ipschitz control-{L}yapunov function.
\newblock In {\em Proceedings of the 39th IEEE Conference on Decision and
  Control}, volume~4, pages 3994--3999. IEEE, 2000.

\bibitem{Khalil1996-nonlin-sys}
H.~Khalil.
\newblock {\em Nonlinear {S}ystems}.
\newblock Prentice-Hall. 2nd edition, 1996.

\bibitem{Ledyaev1997-stabilization-meas-err}
Y.~S. Ledyaev and E.~Sontag.
\newblock A remark on robust stabilization of general asymptotically
  controllable systems.
\newblock In {\em Proc. of Conf. on Information Sciences and Systems, Johns
  Hopkins, Baltimore}, volume 246, page 251, 1997.

\bibitem{Sontag1990-stabilization-survey}
E.~Sontag.
\newblock Feedback stabilization of nonlinear systems.
\newblock In {\em Robust control of linear systems and nonlinear control},
  pages 61--81. Springer, 1990.

\bibitem{Sontag1999-stabilization-disturb}
E.~Sontag.
\newblock Stability and stabilization: discontinuities and the effect of
  disturbances.
\newblock In {\em Nonlinear {A}nalysis, {D}ifferential {E}quations and
  {C}ontrol}, pages 551--598. Springer, 1999.

\bibitem{Sontag1995-nonsmooth-CLF}
E.~Sontag and H.~Sussmann.
\newblock Nonsmooth control-{L}yapunov functions.
\newblock In {\em Proc. of IEEE Conf. on Decision and Control}, volume~3, pages
  2799--2805. IEEE, 1995.

\bibitem{Zalinescu2007-subdiff-cont}
C.~Zalinescu.
\newblock Continuity properties for the subdifferential and
  $\varepsilon$-subdifferential of a convex function and its conjugate.
\newblock {\em Journal of Convex Analysis}, 14(3):479--514, 2007.

\end{thebibliography}

\appendices

\section*{Appendix}

\begin{lemma}
	For all $x \in \ball_R$ for some $R>0$, $0 < \alpha < 1$, and any $\eps>0$, there exists an $\eps$-minimizer $y^{\eps}_{\alpha}(x)$ for \eqref{eqn:inf-conv} satisfying:
	\begin{equation*}
		\| y^{\eps}_{\alpha}(x) - x \| \le \sqrt{2 \bar V} \alpha,
	\end{equation*}
	where $\bar V := \sup_{\|x\| \le R} V(x)$.
	\label{lem:y-eps-near}
\end{lemma}
\begin{IEEEproof}
	Letting $R' := (2 \bar V)^{\nicefrac{1}{2}} \alpha$ yields the following: 
	\begin{align*}
		\inf_{\|x-y\| \le R'} \left( V(y) + \tfrac{1}{2 \alpha^2} \| y-x \|^2 \right)  \le V(x) \le \bar V. 
	\end{align*}
	On the other hand, for any $R'' > R'$, it holds that
	\begin{align*}
		\inf_{R' \le \| x-y \| \le R''} \left( V(y) + \tfrac{1}{2 \alpha^2} \|y-x \|^2 \right) \ge \tfrac{1}{2 \alpha^2} R'^2 \ge \bar V.
	\end{align*}
	Therefore,
	\begin{align*}
		& \inf_{y \in \R^n} \left( V(y) + \tfrac{1}{2\alpha^2} \| y-x \|^2 \right) = \\
		& \inf_{\|x-y\| \le R'} \left( V(y) + \tfrac{1}{2 \alpha^2} \| y-x \|^2 \right)
	\end{align*}
	and the conclusion follows.
\end{IEEEproof}

\begin{lemma}
	With the conditions of Lemma \ref{lem:y-eps-near}, for any $\eps_1 > 0$, an $\eps > 0$ and $0 < \alpha < 1$ for the approximate minimizers $y^{\eps}_{\alpha}(x)$ can be chosen so as to satisfy:
	\begin{equation*}
		V_{\alpha}(x) \le V(x) \le V_{\alpha}(x) + \eps_1.
	\end{equation*}
	\label{lem:V-between}
\end{lemma}
\begin{IEEEproof}
	The first inequality $V_{\alpha}(x) \le V(x)$ follows directly from \eqref{eqn:inf-conv}. As for the second, Lemma \ref{lem:y-eps-near} implies $\| y^{\eps}_{\alpha}(x) - x \| \le (2 \bar V)^{\nicefrac{1}{2}} \alpha$. Let $(2 \bar V)^{\nicefrac{1}{2}} \alpha \le \omega_V \left( \tfrac{\eps_1}{2} \right)$. Choose $\eps < \tfrac{\eps_1}{2}$. It follows from \eqref{eqn:y-eps-def}, that $V(y^{\eps}_{\alpha}) \le V_{\alpha}(x) \ge  + \eps$. By the continuity of $V$, it holds that $V(x) \le V(y^{\eps}_{\alpha}) + \tfrac{\eps_1}{2}$. Putting this together gives the result. 
\end{IEEEproof}

\end{document}